\documentstyle[12pt]{article}

\begin{document}

\title{ Difference between three quantities}
\author{
Robert M. Yamaleev\\
Joint Institute for Nuclear Research, LIT, Dubna, Russia.\\
Universidad Nacional Autonoma de Mexico, Mexico.\\
Email: yamaleev@jinr.ru } \maketitle
%
%
\begin{abstract}
The notion of difference between three and more quantities is
introduced. The method is based on one of the remarkable
properties of the Vandermonde's determinant.

\end{abstract}


{\bf Introduction}.

The notion {\it difference between two quantities} $a$ and $b$
given by $(a-b)$ plays a basic role in mathematics, consequently
in all branches of human activity where the mathematics is
applied. However the long stand question is:  what is {\it the
difference between three (or more) quantities}? This question
frequently arises, for example, in the physics  of systems
consisting of many particles, in economics etc.

 The binary operation
$[a,b]=(a-b)$ possesses the following principal feature: with
respect to the third quantity $c$ this operation is decomposed
into a sum of the same operations between $a$ and $c$, and $c$ and
$b$, i.e.,
$$
[a,b]=[a,c]+[c,b].
$$

There were several problems of mathematics and physics where
investigators needed in the notion of {\it the difference between
three quantities}. Y.Nambu \cite{Nambu} quantizing the generalized
Poisson structure on three dimensional phase space met the problem
of extension of the notion of commutator. In fact,the notion of
commutator can be defined only for the pair of operators $A$ and
$B$ as a difference $AB$ and $BA$. Thus, if one wants to extend
this notion for triple operators he will need on the notion of the
difference between three quantities. In
Refs.\cite{Yamal1},\cite{Yamal2}, the following approach has been
developed. By noting that in formula $(a-b)=(a+\theta b$ the value
$\theta=-1$ is a primitive root of quadratic polynomial $x^2-1$,
the author suggested the following definition of the difference
between three quantities
$$
[a,b,c]=a+b\theta+c\theta^2,
$$
where $\theta$ is a primitive root of polynomial $x^3-1$. In
Refs.\cite{Kerner},\cite{Himbert}, this formula of difference have
been used in order to formulate a notion of ternary commutator.
Apparently, this definition belong to the field of complex numbers
and possesses with the following feature
$$
[a,b,c]=[a,d,f]+[d,a,f]+[d,f,c].
$$

In the present paper we suggest a definition of the notion of the
difference between three and more quantities making use of a
feature of the Vandermonde determinant. Denote by $[a,b,c]$
difference between three quantities $a,b,c$. With respect to
additional quantity $d$ this definition of the difference is
decomposed as follows
$$
[a,b,c]=[d,b,c]+[a,d,c]+[a,b,d].
$$
All quantities belong to the field of real numbers.

\section{ Difference between three quantities }

Let us start with the fraction of type
$$
\frac{1}{x^3-3p_1x^2+2p_2x-p^2}=\frac{1}{(x-x_3)(x-x_2)(x-x_1)},
\eqno(1.1)
$$
where $x_i,i=1,2,3$ are roots of the cubic polynomial
$$
P_3(x):=x^3-3p_1x^2+2p_2x-p^2.  \eqno(1.2)
$$
The following expansion for that fraction holds true
$$
\frac{1}{x^3-3p_1x^2+2p_2x-p^2}=
\frac{(x_3-x_2)}{V}\frac{1}{x-x_1}+\frac{(x_1-x_3)}{V}\frac{1}{x-x_2}+\frac{(x_2-x_1)}{V}\frac{1}{x-x_3},
\eqno(1.3)
$$
where by $V$ we denoted the Vandermonde's determinant of matrix of
order $(3\times 3)$
$$
V=(x_1-x_2)(x_2-x_3)(x_3-x_1)=Det~ \left( \begin{array}{ccc}
1&1&1\\
x_1&x_2&x_3\\
x_1^2&x_2^2&x_3^2
\end{array} \right). \eqno(1.4)
$$

Collecting the terms  in the right-hand side of (1.3) we obtain
$$
\frac{(x_3-x_2)}{V}\frac{1}{x-x_1}+\frac{(x_1-x_3)}{V}\frac{1}{x-x_2}+\frac{(x_2-x_1)}{V}\frac{1}{x-x_3}=
$$
$$
\frac{1}{V}\frac{1}{(x-x_3)(x-x_2)(x-x_1)}(~
$$
$$
(x_3-x_2)(x-x_2)(x-x_3)+(x-x_2)(x-x_1)(x_1-x_2)+(x-x_3)(x-x_1)(x_3-x_1)~).\eqno(1.5)
$$
From this equality we come to the conclusion, that
$$
(x_3-x_2)(x-x_2)(x-x_3)+(x-x_2)(x-x_1)(x_1-x_2)+(x-x_3)(x-x_1)(x_3-x_1)=V.\eqno(1.6)
$$
This equation displays an interesting feature of Vandermonde's
determinant, which to our knowledge still has not been revealed
\cite{Vein}.

We suggest to use formula (1.6) as a formula of difference between
three amounts. The interval of this formal distance is bounded by
three points $x_1,x_2,x_3$ and it is divided into three parts by
using only one point, $x$, among them. Thus, the following formula
$$
V=(x_1-x_2)(x_2-x_3)(x_3-x_1)~ \mbox{ is an analogue of the
interval}~(x_1-x_2),$$ and the property given by formula (1.6) is
an analogue of the following property of the interval between two
points
$$
x_1-x_2=(x_1-x)+(x-x_2).\eqno(1.7)
$$

{\bf The Proof of formula (1.6)}.\\

Formula (1.6) is a consequence of one of the properties of
Vandermonde's   determinant. We have to prove that
$$
Det~ \left( \begin{array}{ccc}
1&1&1\\
x_1&x_2&x_3\\
x_1^2&x_2^2&x_3^2
\end{array} \right)=
$$
$$
Det \left( \begin{array}{ccc}
1&1&1\\
x&x_2&x_3\\
x^2&x_2^2&x_3^2
\end{array} \right)+
Det \left( \begin{array}{ccc}
1&1&1\\
x_1&x&x_3\\
x_1^2&x^2&x_3^2
\end{array} \right)+
Det\left( \begin{array}{ccc}
1&1&1\\
x_1&x_2&x\\
x_1^2&x_2^2&x^2
\end{array} \right).\eqno(1.8)
$$

Consider the following matrix
$$
AV:= \left( \begin{array}{ccc}
1+1&1+1&1+1\\
x_1+x&x_2+x&x_3+x\\
x_1^2+x^2&x_2^2+x^2&x_3^2+x^2
\end{array} \right), \eqno(1.9)
$$
and calculate determinant of this matrix in two ways.

Firstly, let us calculate the determinant on making use of the
method of expansion with respect to lines of the matrix. In this
way we find that
$$
Det(AV) =Det~ \left( \begin{array}{ccc}
1+1&1+1&1+1\\
x_1+x&x_2+x&x_3+x\\
x_1^2+x^2&x_2^2+x^2&x_3^2+x^2
\end{array} \right)
$$
$$
=Det~ \left( \begin{array}{ccc}
1+1&1+1&1+1\\
x_1&x_2&x_3\\
x_1^2+x^2&x_2^2+x^2&x_3^2+x^2
\end{array} \right)+
Det~ \left( \begin{array}{ccc}
1+1&1+1&1+1\\
x&x&x\\
x_1^2+x^2&x_2^2+x^2&x_3^2+x^2
\end{array} \right)
$$
$$
=Det~ \left( \begin{array}{ccc}
2&2&2\\
x_1&x_2&x_3\\
x_1^2+x^2&x_2^2+x^2&x_3^2+x^2
\end{array} \right)=
Det~ \left( \begin{array}{ccc}
2&2&2\\
x_1&x_2&x_3\\
x_1^2&x_2^2&x_3^2
\end{array} \right)+
Det~ \left( \begin{array}{ccc}
2&2&2\\
x_1&x_2&x_3\\
x^2&x^2&x^2
\end{array} \right).
$$
The last determinant is equal to zero. In this way we get
$$
Det(AV)=2V. \eqno(1.20)
$$

Secondly, let us expand the determinant with respect to columns.

For that purpose it is convenient to use the following notation of
the Vandermonde determinant:
$$
V=Det~ \left( \begin{array}{ccc}
1&1&1\\
x_1&x_2&x_3\\
x_1^2&x_2^2&x_3^2
\end{array} \right)=[x_1,x_2,x_3].
$$
In this notation $Det(AV)=[x_1+x,x_2+x,x_3+x]$. Then,  the
following expansion holds true
$$
[x_1+x,x_2+x,x_3+x]=[x_1,x_2+x,x_3+x]+[x,x_2+x,x_3+x]=
$$
$$
[x_1,x_2,x_3+x]+[x_1,x,x_3+x]+[x,x_2,x_3]=
[x_1,x_2,x_3]+[x_1,x_2,x]+[x_1,x,x_3]+[x,x_2,x_3]. \eqno(1.21)
$$
Notice, the first term is the Vandermond's determinant. Therefore,
$$
Det(AV)=V+[x_1,x_2,x]+[x_1,x,x_3]+[x,x_2,x_3]=2V.
$$
Hence,
$$
V=[x_1,x_2,x]+[x_1,x,x_3]+[x,x_2,x_3]. \eqno(1.22)
$$

{\bf End of proof.}

\section{ Difference between $n\geq 2$ quantities}

Since we have found the concept of difference between three
quantities this generalization to the case of $n\geq 3$ quantities
is straightforward.

Consider $n$-th order Vandermonde's matrix
$$
V_{ik}:=\left( \begin{array}{ccccc}
1&1&1&...&1\\
x_1&x_2&x_3&...&x_n\\
x_1^2&x_2^2&x_3^2&...&x_n^2\\
x_1^3&x_2^3&x_3^3&...&x_n^3\\
...&...&...&...&...\\
x_1^{n-1}&x_2^{n-1}&x_3^{n-1}&...&x_n^{n-1}
\end{array} \right). \eqno(2.1)
$$
The determinant of this matrix is given by well-known
Vandermonde's formula:
$$
V=Det(V_{ij})=\prod_{i>k}(x_i-x_k). \eqno(2.2)
$$

Consider the following auxiliary matrix
$$
AV(x):=\left( \begin{array}{ccccc}
1+1&...&1+1&...&1+1\\
x_1+x&...&x_k+x&...&x_n+x\\
x_1^2+x^2&...&x_k^2+x^2&...&x_n^2+x^2\\
...&...&...&...&...\\
x_1^l+x^l&...&x_k^l+x^l&...&x_n^l+x^l\\
...&...&...&...&...\\
x_{n-1}+x^{n-1}&...&x_k^{n-1}+x^{n-1}&...&x_n^{n-1}+x^{n-1}
\end{array} \right).  \eqno(2.3)
$$
Now, let us prove that the determinant of this matrix equal to
$2V$. Firstly,  expand this determinant with respect to $l$-th
line. The result is given by sum of two determinants
$$
Det(AV(x))=Det(AV_1(x,x_k^l))+Det(AV_2(x,x^l)), \eqno(2.4)
$$
where we denoted
$$
AV_1(x,x^l)=\left( \begin{array}{ccccc}
1+1&...&1+1&...&1+1\\
x_1+x&...&x_k+x&...&x_n+x\\
x_1^2+x^2&...&x_k^2+x^2&...&x_n^2+x^2\\
...&...&...&...&...\\
x_1^l&...&x_k^l&...&x_n^l\\
...&...&...&...&...\\
x_1^{n-1}+x^{n-1}&...&x_k^{n-1}+x^{n-1}&...&x_n^{n-1}+x^{n-1}
\end{array} \right),
$$
and,
$$
AV_2(x,x^l)=\left( \begin{array}{ccccc}
1+1&...&1+1&...&1+1\\
x_1+x&...&x_k+x&...&x_n+x\\
x_1^2+x^2&...&x_k^2+x^2&...&x_n^2+x^2\\
...&...&...&...&...\\
x^l&...&x^l&...&x^l\\
...&...&...&...&...\\
x_1^{n-1}+x^{n-1}&...&x_k^{n-1}+x^{n-1}&...&x_n^{n-1}+x^{n-1}
\end{array} \right).
$$
The determinant of the second matrix is trivial because there
$l$-th line is proportional to the first one. Continue to expand
the first determinant $AV_1(x,x^l)$ with respect to other lines.
At the final step of this process the determinant of the auxiliary
matrix is reduced to the following form
$$
Det(AV(x))=Det\left( \begin{array}{ccccc}
1+1&...&1+1&...&1+1\\
x_1&...&x_k&...&x_n\\
x_1^2&...&x_k^2&...&x_n^2\\
...&...&...&...&...\\
x_1^l&...&x_k^l&...&x_n^l\\
...&...&...&...&...\\
x_1^{n-1}&...&x_k^{n-1}&...&x_n^{n-1}
\end{array} \right),  \eqno(2.5)
$$
which obviously equal to $2V$,
$$
Det(AV(x))=2V. \eqno(2.6)
$$

Now let us calculate the  determinant $Det(AV(x))$ by expanding
with respect to columns. For the sake of convenience denote the
Vandermonde's determinant (2.1) as follows
$$
V=Det(V[ij])=[x_1...x_k...x_n]. \eqno(2.7)
$$
Correspondingly, the determinant of the auxiliary matrix will be
written in the form
$$
Det(AV(x))=[x_1+x...x_k+x...x_n+x]. \eqno(2.8)
$$
The expansion process with resect to columns of this determinant
is worked out as follows.
$$
[x_1+x...x_k+x...x_n+x]=[x_1,x_2+x...x_k+x...x_n+x]+[x,x_2+x...x_k+x...x_n+x].
\eqno(2.9)
$$
The second term in right-hand side is equal to
$$
[x,x_2+x...x_k+x...x_n+x]=[x,x_2...x_k...x_n].
$$
Continue to expand the first of the sum,
$$
[x_1,x_2+x...x_k+x...x_n+x]=[x_1,x_2,x_3+x...x_k+x...x_n+x]+[x,x_2,x...x_k+x...x_n+x].
$$
The last term is equal to zero. The first one is represented as
follows
$$
[x_1,x_2,x,x_4+x...x_k+x...x_n+x]=[x,x_2,x_3,x_4+x...x_k...x_n].
$$
At the final step of this process  we come to the following
equation
$$
Det(AV(x))=2V=\sum^n_{k=1}[x_1,x_2,...x_{k-1},x,x_{k+1}...x_n]+V,
\eqno(2.10)
$$
On the other hand, according to (2.6) $Det(AV(x))=2V$. Hence,
$$
V=\sum^n_{k=1}[x_1,x_2,...x_{k-1},x,x_{k+1}...x_n]. \eqno(2.11)
$$

This formula implies one of the important features of the
Vandermond's determinant. That is the formula which we suggest to
use as a definition of the {\it difference between $n$
quantities}.

{\bf Concluding remarks}.

Ternary algebraic operations and cubic relations have been
considered, although quite sporadically, by several authors
already in the XIX-th century, e.g. by A. Cayley (\cite{Cayley})
and J.J. Sylvester ( \cite{Sylvester}. The development of Cayley's
ideas, which contained a cubic generalization of matrices and
their determinants, can be found in a recent book by M. Kapranov,
I.M. Gelfand and A. Zelevinskii (\cite{Kapranov}). A discussion of
the next step in generality, the so called $n-ary$ algebras, can
be found in (\cite{Vainerman}).

The difference between two quantities has direct geometrical
interpretation as a distance between two points on a straight
line. Let $O,A,B$ be a set of points on the line and let point $O$
be a point on the left-hand side of the points $A$ and $B$. Let
the values $d(OA),d(OB)$ mean distances between points $A$ and $B$
of the point $O$, correspondingly. Then the difference
$[d(OA),d(OB)]$ does not depend of the motion of the point $O$ and
means the distance between points $A$ and $B$.

In the similar way, let $O$ be a point on the straight line on the
left-hand side of three points $A,B,C$ installed on the same line.
Let $d(OA),d(OB),d(OC)$ be distances from $O$ till points $A,B,C$,
correspondingly. Then the difference $[d(OA),d(OB),d(OC)]$ does
not depend of the motion of the point $O$ along the straight line.
This definition of the difference we suggest use in geometry in
order to found a concept of ternary distance between three points.
A generalization to the case of $n\geq 3$ points is
straightforward.

{99}


\begin{thebibliography}{99}

\bibitem{Nambu} Nambu Y., Generalized Hamilton dynamics. Physical Review D {\bf 7}, p.2405 (1973)
\bibitem{Yamal1} Yamaleev R.M. On geometrical models in 3D space
with cubic metrics. Communications of Joint Institute for Nuclear
Research, P5-89-269, Dubna, 1989.
\bibitem{Yamal2} Yamaleev R.M. Fractional power of momenta and
paragrassmann extension of Pauli equation. Adv.Appl.Clifford
Alg.{\bf 7} (S) (1997)279.
\bibitem{Himbert} A. Himbert, Comptes Rendus de l'Acad.Sci. Paris,
(1935).
\bibitem{Kerner} R. Kerner, {\it ``The Cubic Chessboard''}, Class. and Quantum
Gravity, {\bf 14} 1A, p. A203 (1997).
\bibitem{Vein} Vein R., Dale P., "Determinants and their
applications in mathematical physics", Springer-Verlag, New York,
Inc., 1999. ISBN 0-387-98558-1.
\bibitem{Cayley} A. Cayley , Cambridge Math. Journ. {\bf 4}, p. 1 (1845)
\bibitem{Sylvester} J.J. Sylvester, Johns Hopkins Circ. Journ., {\bf 3},
p.7 (1883).
\bibitem{Kapranov} M. Kapranov, I.M. Gelfand, A. Zelevinskii, {\it Discriminants,
Resultants, and Multidimensional Determinants}, Birkh\"auser ed.,
(1994)
\bibitem{Vainerman} L. Vainerman, R. Kerner, Journal of Math. Physics,
{\bf 37} (5), p. 2553 (1996)




\end{thebibliography}
\end{document}